\documentstyle[12pt,amscd]{amsart}
\newtheorem{Theorem}{Theorem}[section]
\newtheorem{Lemma}[Theorem]{Lemma}
\newtheorem{Corollary}[Theorem]{Corollary}

\def\depth{\operatorname{depth}}
\def\reg{\operatorname{reg}}
\def\geom{\operatorname{g-reg}}
\def\reltype{\operatorname{reltype}}
\def\To {\longrightarrow}

\def\mm{{\frak m}}
\def\pp{{\frak p}}

\begin{document}
\title{Uniform bounds in\\ generalized Cohen-Macaulay rings}
\author{Cao Huy Linh and Ngo Viet Trung}
\address{Department of  Mathematics,  University of Education, 32 Le Loi, Hue, Vietnam}
\email{huylinh2002@@yahoo.com}
\address{Institute of  Mathematics,  18 Hoang Qu\^oc Vi\^et, 10307 Hanoi, Vietnam}
\email{nvtrung@@math.ac.vn}

\date{}
\maketitle
\begin{abstract}
We establish a uniform bound for the Castelnuovo-Mumford
regularity of associated graded rings of parameter ideals in a generalized Cohen-Macaulay ring. 
As consequences, we obtain uniform bounds for the relation type and the postulation number. 
Moreover, we show that generalized Cohen-Macaulay rings can be characterized by the existence of such uniform bounds.  
\end{abstract}

\section*{Introduction} 

Let $R$ be a finitely generated standard graded algebra over a local ring. Let $R_+$ be the ideal of $R$ generated by the elements of positive degrees of $R$. The Castelnuovo-Mumford regularity of $R$ is the number
$$\reg(R) := \max\{a(H_{R_+}^i(R))+i|\ i \ge 0\},$$
where $a(H_{R_+}^i(R))$ denotes the largest non-vanishing degree of the local cohomology module $H_{R_+}^i(R)$ of $R$ with respect to $R_+$.\smallskip

For any ideal $I$ in a local ring $A$ let $G_I(A)$ denote the associated graded ring of $A$ with respect to $I$. One can use the Castelnuovo-Mumford regularity $\reg(G_I(A))$  to study the behavior of ideal powers of $I$. In fact,  the reduction number $r(I)$ of $I$ is bounded by $\reg(G_I(A))$  and the relation type $\reltype(I)$ of $I$ is bounded by $\reg(G_I(A))+1$. 
More precisely, $\reg(G_I(A))$ is the maximum of a finite number of invariants of $I$ which include $r(I)$ and $\reltype(I)-1$ \cite{Tru}, \cite{Tru2}. Moreover, if $I$ is an ideal of finite colength,
the postulation number of $I$ is bounded by $\reg(G_I(A))+1$. \smallskip

Recall that a local ring $(A,\mm)$ is called  {\it generalized Cohen-Macaulay} if the local cohomology module $H_\mm^i(A)$ has finite length for $i < \dim A$ \cite{CST}.
This class of rings is rather large because it contains the local rings at the vertices of affine cones over projective non-singular varieties and the local rings of isolated singularities. \smallskip

The aim of this paper is to prove that there exists a uniform  bound for  $\reg(G_Q(A))$ for all parameter ideals $Q$ of a generalized Cohen-Macaulay ring $A$. 
This implies a uniform bound for the relation type  of $Q$ which was  established first by Lai \cite{La} in the case the residue field of $A$ is finite and then by Wang \cite{W2} for any generalized Cohen-Macaulay local ring. This also implies a uniform bound for the postulation number of $Q$ which was known only in the case $\dim A = 2$ \cite{W1}. It should be mentioned that uniform bounds for the relation type are closely related to uniform bounds for the Artin-Rees number which 
plays an important role in several topics of commutative algebra \cite{DO}, \cite{H1}, \cite{H2}, \cite{OC}, \cite{P}. \smallskip

One can associated with a generalized Cohen-Macaulay local ring $A$ the invariant
$$I(A) = \sum_{i=0}^{d-1}  {d-1 \choose i} \ell(H_\mm^i(A))$$
where $d = \dim A$. If $d > 0$,  we prove that
\begin{align*}
\reg(G_Q(A)) & \le \max\{I(A)-1,0\}\ \text{if}\ d = 1,\\
\reg(G_Q(A)) & \le  \max\{(4I(A))^{(d-1)!}-I(A)-1,0\}\ \text{if}\ d \ge 2.
\end{align*}
for all parameter ideals $Q$ of $A$. From these bounds we immediately obtain
\begin{align*}
\reltype(Q) & \le \max\{I(A),1\}\ \text{if}\ d = 1,\\
\reltype(Q)  & \le  \max\{(4I(A))^{(d-1)!}-I(A),1\}\ \text{if}\ d \ge 2.
\end{align*}
Moreover, the same bounds also hold if we replace $\reltype(Q)$ by the postulation number of $Q$.  \smallskip

Bounds for the Castelnuovo-Mumford regularity of associated graded rings of $\mm$-primary ideals were given already in \cite{Li}, \cite{RTV1}, \cite{Tri} but in terms of the (extended) degree, which depends on the given ideal and hence are not uniform. Our approach follows the method of \cite{RTV1} which is based on a result of Mumford which links  the geometric regularity to that of a  general hyperplane section by means of the difference of the Hilbert polynomial and the Hilbert function. It turns out that this difference can be uniformly bounded for the associated graded rings of parameter ideals in a generalized Cohen-Macaulay ring. 
The key point is to show that there is a uniform bound for the invariant $I(A/J^n)$ for all ideals $J$ generated by subsystems of parameters of $A$. 
\smallskip

We will also discuss the problem which local rings have uniform bounds for the relation type of parameter ideals. This problem originated from Huneke's question whether equidimensional unmixed local rings have the above property. Recently, Aberbach found a counter-example  (see \cite{AGH}). Results of Wang \cite{W1} and  of Aberbach, Ghezzi and Ha \cite{AGH} showed that the class of rings with the above property must be larger than the class of generalized Cohen-Macaulay rings. However, we will prove that a local ring is generalized Cohen-Macaulay if and only if 
there exists an integer $r$ such that for  every quotient ring by an ideal generated by a subsystem of parameters, the relation type (or the regularity of the associated graded rings) of parameter ideals is bounded by $r$. Therefore, if there exists a uniform bound for the relation type of parameter ideals, then the bound must depend on invariants which may increase when we pass to quotient rings of
ideals generated by subsystems of parameters. That helps explain why it is so difficult to characterize local rings with the above property. 

\section{Generalized Cohen-Macaulay rings}

The aim of this section is to establish properties of generalized Cohen-Macaulay rings which will be used in the proofs of the main results of this paper. \smallskip

Let $(A,\mm)$ be a local ring and $d = \dim A$. 
For every parameter ideal $Q$ of $A$ we define
$$I(Q,A) := \ell(A/Q) - e(Q,A),$$
where $e(Q,A)$ denotes the multiplicity of $A$ with respect to $Q$.
It is known that $A$ is a generalized Cohen-Macaulay ring if and only if
$\sup I(Q,A) < \infty$ \cite[(3.3)]{CST}. We always have
$$I(A) = \sup I(Q,A).$$
In particular, for a fixed parameter ideal $Q = (x_1,...,x_d)$, the integers $I((x_1^n,...,x_d^n),A)$ form a non-decreasing sequence and
$I(A) = \sup I((x_1^n,...,x_d^n),A)$ \cite[Theorem 2.1 and Lemma 2.2]{Tru1}.
This provides an effective way for the computation of $I(A)$. \smallskip

If $A$ is a generalized-Cohen-Macaulay ring, then 
$\dim A/\pp = d-i$ for all associated prime $\pp \neq \mm$ of $(x_1,...,x_i)$, $i = 1,...,d-1$
\cite[(2.5)]{CST}. By a result of Auslander and Buchsbaum \cite[Corollary 4.8]{AB}, this implies 
$$e(Q,A/(x_1,...,x_i)) = \ell(A/Q)  - \ell((x_1,...,x_{d-1}):x_1/(x_1,...,x_{d-1})) = e(Q,A).$$
Hence $I(Q, A/(x_1,...,x_i)) = I(Q,A)$. Furthermore, $A/(x_1,...,x_i))$ is a generalized Cohen-Macaulay ring with 
$$I(A/(x_1,...,x_i)) \le I(A).$$

Now we shall see that the Hilbert function of $Q$ can be bounded solely in terms of $e(Q,A)$ and $I(A)$.

\begin{Lemma}\label{Hilbert} 
Let $A$ be a generalized Cohen-Macaulay ring with $\dim A = d > 0$. Let $Q$ be an arbitrary parameter ideal of $A$. For all $n \ge 0$,
$$\ell(A/Q^{n+1}) \le \binom{n+d}{d}e(Q,A)+\binom{n+d-1}{d-1}I(A). $$
\end{Lemma}

\begin{pf}
If $d = 1$, we use the inequality 
$$\ell(A/Q^{n+1}) \le  \ell(A/Q^{n+1}+L) + \ell(L),$$
where $L$ denotes the largest ideal of finite length of $A$. Since $A/L$ is an one-dimensional Cohen-Macaulay ring and $L = H_\mm^0(A)$, 
$$\ell(A/Q^{n+1}+L) + \ell(L) = (n+1)e(Q,A) + I(A),$$
which proves this case.\par
If $d > 1$, we put $x = x_1$. Then  $e(Q,A/(x)) = e(Q,A)$  and $I(A/(x)) \le
I(A)$. Using induction on $d$ we may assume that for all $i \ge 0$,
$$\ell(A/(Q^{i+1},x)) \le {i+d-1 \choose d-1}e(Q,A) + {i+d-2
\choose d-2}I(A).$$
From the exact sequence
$$0 \To Q^{i+1}:x/Q^i  \To A/Q^i \overset {x}  \To
A/Q^{i+1} \To A/(Q^{i+1},x) \To 0$$
we deduce that
$$\ell(Q^i/Q^{i+1}) = \ell(A/Q^{i+1})-\ell(A/Q^i) \le
\ell(A/(Q^{i+1},x)).$$
Using the above inequalities we get
\begin{align*}
\ell(A/Q^{n+1}) & = \sum_{i=0}^n \ell(Q^i/Q^{i+1}) \le
\sum_{i=0}^n \ell(A/(Q^{i+1},x))\\
& \le \sum_{i=0}^n \left[{i+d-1 \choose d-1}e(A) + {i+d-2 \choose
d-2}I(A)\right]\\
& = {n+d \choose d}e(A) + {n+d-1 \choose d-1}I(A).
\end{align*}
\end{pf}

\noindent{\it Remark.} There is a more precise upper bound for $\ell(A/Q^{n+1})$  in \cite[Theorem 4.1]{Tru1}.\medskip

Using the above lemma we can show that there exists a uniform bound for the invariants $I(A/J^{n+1})$, where $J$ is any ideal generated by a subsystem of parameters of $A$.

\begin{Theorem} \label{invariant}
Let $A$ be a generalized Cohen-Macaulay ring with $\dim(A) = d$. Let $x_1,...,x_i$ be a subsystem of parameters of
$A$ and $J = (x_1,...,x_i)$, $0< i < d$. Then $A/J^{n+1}$ is a generalized Cohen-Macaulay ring with
 $$I(A/J^{n+1}) \leq  \binom{n+i-1}{i-1}I(A).$$
\end{Theorem}

\begin{pf}
We extend $x_1,...,x_i$ to a system of parameters $x_1,...,x_d$ of $A$ and define $B := A/(x_{i+1},...,x_d)$. Then
$$\ell(A/(J^{n+1},x_{i+1},...,x_d)) = \ell(B/J^{n+1}B).$$
By  Lemma \ref{Hilbert},
$$
\ell(B/J^{n+1}B)  \le \binom{n+i}{i}e(J,B) +
\binom{n+i-1}{i-1}I(B).
$$
Put $Q = (x_1,...,x_d)$. Then $e(J,B) = e(Q,A)$ (here we need $0< i < d$) and $I(B) \le I(A)$. Therefore,
$$
\ell(A/(J^{n+1},x_{i+1},...,x_d))  \le \binom{n+i}{i}e(Q,A) +
\binom{n+i-1}{i-1}I(A).
$$

By the associative formula for multiplicity, 
$$ 
e((x_{i+1},...,x_d),A/J^{n+1}) =
\sum_{\pp}\ell(A_{\pp}/J^{n+1}A_{\pp})e((x_{i+1},...,x_d),A/\pp),
$$
  where the sum is taken over all prime ideals $\pp \supseteq J$ with
  $\dim(A/\pp) = d-i$. Since $A_{\pp}$ is  a Cohen-Macaulay ring with $\dim A_{\pp} =  i$ and $J A_\pp$ is a parameter  ideal of $A_\pp$ \cite[(4.8)]{CST},
$$
 \ell(A_{\pp}/J^{n+1}A_{\pp}) = \binom{n+i}{i}\ell(A_{\pp}/J
 A_{\pp}).
$$
Therefore,
\begin{align*} 
e((x_{i+1},...,x_d),A/J^{n+1})  & =
 \binom{n+i}{i}\sum_{\pp}\ell(A_{\pp}/J A_{\pp})e((x_{i+1},...,x_d),A/\pp)\\
&  = \binom{n+i}{i}e((x_{i+1},...,x_d) ,A/J)\\
& = \binom{n+i}{i}e(Q,A).
 \end{align*}

It follows that
\begin{align*}
I((x_{i+1},...,x_d),A/J^{n+1})  & =   \ell(A/(J^{n+1},x_{i+1},...,x_d)) -
e((x_{i+1},...,x_d),A/J^n)\\
& \leq \binom{n+i-1}{i-1}I(A).
\end{align*}
This implies 
$$I(A/J^{n+1}) = \sup I((x_{i+1},...,x_d),A/J^{n+1}) \leq \binom{n+i-1}{i-1}I(A).$$
\end{pf}

\begin{Corollary}\label{Wang1}
Let $A$ be a generalized Cohen-Macaulay ring with $\dim(A) = d$. Let $x_1,...,x_{i+1}$ be a subsystem of parameters of
$A$ and $J = (x_1,...,x_i)$, $0 < i < d$. For all $n \ge 0$ and $m \ge 1$,
$$\ell(J^{n+1}:x_{i+1}^m/J^{n+1}) \leq  \binom{n+i-1}{i-1}I(A).$$
\end{Corollary}

\begin{pf}
Let $\pp \neq \mm$ be any associated prime ideal of $J^{n+1}$, then $\dim A/\pp = \dim A/J^{n+1} = d-i$. Therefore, $x_{i+1} \not\in \pp$. This implies 
$J^{n+1}:x_{i+1}^m/J^{n+1}\subseteq H_\mm^0(A/J^{n+1})$.
Thus,
$$\ell(J^{n+1}:x_{i+1}^m/J^{n+1}) \leq \ell(H_\mm^0(A/J^{n+1})) \le I(A/J^{n+1}).$$
Hence the statement immediately follows from Theorem \ref{invariant}.
\end{pf}

\noindent{\it Remark.} Wang \cite[Theorem 3.3]{W2} already proved that there exists a uniform upper bound for $\ell(J^n:x_{i+1}/J^n)$ in terms of $\ell(H_\mm^i(A))$, $i < d$. However, the bound is not explicit and the proof is very complicated.
\medskip

We shall need  the following result in the next section.

\begin{Corollary} \label{quotient}
Let $A$ be a generalized Cohen-Macaulay ring with $\dim(A) = d
\geq 2$. Let $Q = (x_1,...,x_d)$ be a parameter ideal of
$A$. For all $n \ge 0$ and $m \ge 1$,
$$\ell(Q^{n+m}:x_d^m/Q^n) \le {n+d-2 \choose d-2}I(A).$$
\end{Corollary}

\begin{pf}
Put $J = (x_1,...,x_{d-1})$. Since $Q^{n+m} \subseteq J^{n+1} + x_d^mQ^n$,
\begin{align*}
\ell(Q^{n+m}:x_d^m/Q^n) & \le \ell((J^{n+1}:x_d^m)+Q^n/Q^n)\\
& \le \ell(J^{n+1}:x_d^m/(J^{n+1}:x_d^m) \cap Q^n)\\
& \le \ell(J^{n+1}:x_d^m/J^{n+1})\\
& \le {n+d-2 \choose d-2}I(A).
\end{align*}
\end{pf}

It is known that a local ring $A$ is  generalized Cohen-Macaulay if and only if there exists an integer $n > 0$ such that for every system of parameters
$x_1,...,x_d$ of $A$ and $i = 1,...,d$,
$$(x_1,...,x_{i-1}):x_i \subseteq (x_1,...,x_{i-1}):\mm^n$$
(see \cite[(3.3)]{CST}). We shall  use the following modification of this characterization.

\begin{Lemma} \label{characterization}
$A$ is a generalized Cohen-Macaulay ring if and only if there exists an integer $n > 0$ such that for every system of parameters $x_1,...,x_d$ of $A$,
$$(x_1,...,x_{d-1}):x_d \subseteq (x_1,...,x_{d-1}):\mm^n.$$
\end{Lemma}

\begin{pf}
We only need to show that $A$ is a generalized Cohen-Macaulay ring if there exists such an integer $n$. But this follows from the following simple observation
\begin{align*}
(x_1,...,x_{i-1}):x_i & = \cap_{m \ge 1}(x_1,...,x_{i-1},x_{i+1}^m,...,x_d^m):x_i\\
& \subseteq \cap_{m \ge 1}(x_1,...,x_{i-1},x_{i+1}^m,...,x_d^m):\mm^n\\
& = (x_1,...,x_{i-1}):\mm^n.
\end{align*}
\end{pf}

\section{Uniform bounds for the regularity}

Let $(A,\mm)$ be a local ring and $d = \dim A > 0$.
Let $I$ be an $\mm$-primary ideal of $A$.
First, we will describe a method to estimate the regularity of 
the associated graded ring $G_I(A)$.
This method was used in \cite{RTV1} to prove that there is an upper bound for $\reg(G_\mm(A))$ in terms of the extended degree of $A$.
It was later generalized to give upper bounds for $\reg(G_I(M))$ for an arbitrary finitely generated $A$-module $M$ in \cite{Li}.
\smallskip

Let $L$  denote the largest ideal of finite length of $A$. By \cite[Lemma 4.3]{Li} or \cite[Lemma 3.1]{RTV1} we have
$$\reg(G_I(A)) \le \reg(G_I(A/L)) + \ell(L). \eqno{(1)}$$
Therefore, we only need to estimate $\reg(G_I(A/L))$.
Since $\depth(A/L) > 0$, 
we may assume from the beginning that $\depth A>0$.
In this case, we have
$$\reg(G_I(A)) =\geom(G_I(A)) \eqno{(2)}$$
by \cite[Lemma 4.2]{Li} or \cite[Proposition 3.2]{RTV1}. \smallskip

Recall that for a finitely generated standard graded algebra $R$ over a local ring,
$$\geom(R) := \max\{a(H_{R_+}^i(R))+i|\ i > 0\}.$$
This invariant is called the {\it geometric regularity} of $R$ because it is related to the regularity of the sheaf associated with $R$ \cite{RTV1}.
It is obvious that $\geom(R) \le \reg(R)$.
Due to a result of Mumford \cite[p.~101, Theorem]{Mu}, we can estimate $\geom(R)$ by means of the geometric regularity of a ``generic hyperplane section".\smallskip

Recall that a homogeneous element $z \in R$ is called {\it filter-regular} if $(0:z)_n = 0$ for $n \gg 0$ \cite{Tru}. This is equivalent to the condition $z \not\in P$ for any associated prime ideal $P \not\supseteq R_+$ of $R$. Therefore, if the local ring $R_0$ has infinite residue field, filter-regular elements of any degree exist.

\begin{Theorem}\label{Mumford} {\rm (see \cite[Theorem 2.7]{Li} or \cite[Theorem 1.4]{RTV1})}
Let $R$  be a finitely generated standard graded algebra over an artinian local ring.
Let  $z \in R_1$ be a filter-regular element.
Let $n$ be an integer such that $ \geom(R/(z)) \le n$. Then  
$$\geom(R) \le n + p_R(n) - h_{R/L}(n),$$
where $p_R(n)$ is the Hilbert polynomial of $R$, $L$  is the largest ideal of finite length of $R$ and $h_{R/L}(n)$ is the Hilbert function of $R/L$.
\end{Theorem}

By a result of Serre we always have
$$h_R(n) - p_R(n) = \sum_{i=0}^{\dim R}\ell(H_{R_+}^i(R)_n)$$
(see e.g. \cite[Theorem 17.1.6]{BS}). From this it follows that
$h_R(n) = p_R(n)$ for all $n \ge \reg(R)+1$.\medskip

The least number $n_0$ such that $h_{G_I(A)}(n) =  p_{G_I(A)}(n)$ for $n \ge n_0$ is called the {\it postulation number} of $I$, denoted by $p(I)$. As we have seen above,
$$p(I) \le \reg(G_I(A))+1.$$
In particular, this implies that the Hilbert-Samuel function $\ell(A/I^{n+1})$ is a polynomial for $n \ge  \reg(G_I(A))$.\smallskip

Choose an element $x$ in $I \setminus I^2$ such that its initial form $x^*$ is a filter-regular element in $G_I(A)$. Then 
$$\geom(G_I(A)/(x^*)) = \geom(G_I(A/(x)) \eqno{(3)}$$
by \cite[Lemma 4.1]{Li} or \cite[Lemma 2.2]{RTV1}. Using induction on $\dim A$ we may assume some bound for $ \geom(G_I(A/(x)) \le \reg(G_I(A/(x))$. Thus, to bound $\reg(G_I(A))$  we only need to estimate $p_{G_I(A)}(n)-h_{G_I(A)/L}(n)$ for $n$ greater or equal the assumed bound for $ \geom(G_I(A/(x))$. \smallskip

\begin{Lemma} \label{difference}
Let $A$ and $x$ be as above. Then there exists an integer $m $ such that
$$p_{G_I(A)}(n)-h_{G_I(A)/L}(n) \le \ell(I^{n+1}:x/I^n) + \ell(I^{n+m+1}:x^m/I^{n+1})$$
for all $n \ge \reg(G_I(A/(x))$.
\end{Lemma}

\begin{pf}
We have the exact sequence
$$0 \longrightarrow (I^{n+1}:x)/I^n  \longrightarrow A/I^n\overset{x}
\longrightarrow A/I^{n+1} \longrightarrow 
A/(I^{n+1},x)\longrightarrow 0.$$
From this it follows that
$$h_{G_I(A)}(n) = \ell(I^n/I^{n+1}) =  \ell(A/(I^{n+1},x)) - \ell(I^{n+1}:x/I^n).$$
As observed above, $\ell(A/(I^{n+1},x))$ is a polynomial for $n \ge \reg(G_I(A/(x))$.
For $n \gg 0$ we have $I^{n+1}:x = I^n +(0:x)$ \cite[Lemma 4.4 (i)]{Tru2} and $(0:x) \cap I^n = 0$ \cite[Proposition 4.6]{Tru2}. Therefore, $\ell(I^{n+1}:x/I^n) = \ell((0:x)/(0:x) \cap I^n) = \ell(0:x)$
for $n \gg 0$. For $n \ge \reg(G_I(A/(x))$, this implies 
$$p_{G_I(A)}(n) = \ell(A/(I^{n+1},x)) - \ell(0:x)$$
and hence 
$$p_{G_I(A)}(n) - h_{G_I(A)}(n) \le \ell(I^{n+1}:x/I^n).$$
Since $x^*$ is a filter-regular element of $G_I(A)$, there exists an integer $m$ such that
$L = 0:_{G_I(A)}(x^*)^m$. From this it follows that
$$\ell(L_n) =  \ell((I^{m+n+1}:x^m) \cap I^n /I^{n+1}) \le \ell(I^{m+n+1}:x ^m/I^{n+1}).$$
Therefore,
\begin{align*}
p_{G_I(A)}(n)-h_{G_I(A)/L}(n) & = p_{G_I(A)}(n)-h_{G_I(A)}(n) + \ell(L_n)\\
& \le \ell(I^{n+1}:x/I^n) + \ell(I^{m+n+1}:x^m/I^{n+1}).
\end{align*}
\end{pf}

Now we are ready to establish a uniform bound for the regularity of the associated graded rings.

\begin{Theorem}\label{main}
Let $A$ be a generalized Cohen-Macaulay ring with $d=\dim(A)\geq 1 $. For all parameter ideals $Q$ of $A$, 
\begin{align*}
\reg(G_Q(A)) & \le \max\{I(A)-1,0\}\ \text{if}\ d = 1,\\
\reg(G_Q(A)) & \le \max\{(4I(A))^{(d-1)!}-I(A)-1,0\}\ \text{if}\ d \ge 2.
\end{align*}
\end{Theorem}

\begin{pf}
Let $L$  denote the largest ideal of finite length of $A$. Then $L = H_\mm^0(A)$. \par

If $d = 1$, then $Q = (x)$. Since $L:x = L$, we have $(x^n) \cap L= x^nL$ for all $n \ge 0$. Therefore,
$$G_Q(A/L) = \oplus_{n\ge 0}(x^n,L)/(x^{n+1},L) = \oplus_{n\ge 0}(x^n)/(x^{n+1},x^nL)$$
so that we have an exact sequence of the form
$$0 \to K \to G_Q(A) \to G_Q(A/L) \to 0$$
where $K = \oplus_{n\ge 0}(x^{n+1},x^nL)/(x^{n+1}) = \oplus_{n\ge 0}x^nL/x^{n+1}L.$
It is obvious that $K_n = 0$ for $n \ge \ell(L) = I(A)$. Therefore,
$$\reg(G_Q(A)) \le \max\{I(A)-1,\reg(G_Q(A/L))\}.$$
Since $A/L$ is a Cohen-Macaulay ring, $G_Q(A/L)$ is isomorphic to a polynomial ring over $A/(L,Q)$. Therefore, $\reg(G_Q(A/L)) = 0$, which implies the first bound.\par

If $ d \ge 2$, we first note that
$$I(A) = I(A/L) + \ell(L).$$
Therefore, using (1) we only need to prove the second bound for the local ring $A/L$.
So we may assume that  $\depth A > 0$. By (2) this implies
$$\reg(G_Q(A)) = \geom(G_Q(A)).$$

If $A$ is a Cohen-Macaulay ring, then $G_Q(A)$ is isomorphic to a polynomial ring over $A/Q$, whence $\geom(G_Q(A)) = 0$. Therefore, we may assume that $A$ is not a Cohen-Macaulay ring (i.e. $I(A) \ge 1$) with $\depth A > 0$. \par

Let $Q = (x_1,...,x_d)$ and $x = x_1$. Notice that $A/(x)$ is  not a Cohen-Macaulay ring with $I(A/(x)) \le I(A)$. Without loss of generality we
may further assume that the residue field of $A$ is infinite.  Then we may choose $x$ such that its initial form $x^*$   is a filter-regular element in $G_Q(A)$. By (3) we have
$$\reg(G_Q(A)/(x^*)) = \reg(G_Q(A/(x))).$$
Let $n$ be any integer such that $n \ge \reg(G_Q(A/(x)))$. 
By Theorem \ref{Mumford} and Lemma \ref{difference}, there exists an integer $m$ such that
$$\geom(G_Q(A)) \leq  n+ \ell(Q^{n+1}:x/Q^n) + \ell(Q^{n+m+1}:x^m/Q^{n+1}).$$
Applying Corollary \ref{quotient} we obtain
\begin{align*}
\geom(G_Q(A))  & \le n + {n+d-2 \choose d-2}I(A) + {n+d-1 \choose d-2}I(A)\\
& \le n + (n+1)^{d-2}I(A) + (n+2)^{d-2}I(A)\end{align*}
where the last inequality can be easily checked.\par

If $d = 2$, we may put $n = I(A)-1$ and obtain
$$\geom(G_Q(A)) \le 3I(A)-1 = (4I(A))^{(2-1)!}-I(A)-1.$$

If $d > 2$, using induction on $d$ we may assume that 
$$\geom(G_Q(A/(x))) \le \big(4I(A/(x))\big)^{(d-2)!}-I(A/(x))-1.$$
Since $I(A/x) \le I(A)$, 
$$\big(4I(A/(x))\big)^{(d-2)!}-I(A/(x)) -1\le (4I(A))^{(d-2)!}-I(A)-1.$$
Therefore, we may put 
$n = (4I(A))^{(d-2)!}-I(A)-1$. Note that $n \ge I(A)$ because $I(A) \ge 1$. Then
\begin{align*}
\geom(G_Q(A))  &  \le n + (n+1)^{d-2}I(A) + (n+2)^{d-2}I(A)\\
& \le nI(A) + (n+1)^{d-2}I(A) + \sum_{i=0}^{d-2}{d-2 \choose i}(n+1)^{d-2-i}I(A)\\
& \le  \sum_{i=0}^{d-2}{d-1 \choose i+1}(n+1)^{d-2-i}I(A)\ (\text{because}\ d > 2)\\
& \le  \sum_{i=1}^{d-1}{d-1 \choose i}(n+1)^{d-1-i}I(A)^i \\
& \le (n+1 + I(A))^{d-1}- I(A)-1\ (\text{because}\ n \ge I(A)) \\
& \le  (4I(A))^{(d-1)!}-I(A)-1.
\end{align*}
\end{pf}

The following example shows that  the bound for $\reg(G_Q(A))$ in the case $d=1$ is sharp.  \medskip

\noindent{\it Example.}
Let $A = k[[x,y]]/(x^2,xy^r)$, $r \ge 1$, and $Q = yA$. Then
$I(A) = \ell(H_\mm^0(A)) = \ell(xA) = r$. It is easy to check that 
$$G_Q(A) \cong A[T]/(y,xT^r).$$
Since the maximal degree of the generators of $G_Q(A)$ is bounded by $\reg(G_Q(A))+1$,
$r-1 \le \reg(G_Q(A))$. By Theorem \ref{main}, this implies $\reg(G_Q(A)) = r-1$.
\medskip

As one can see from the above proof,  the bound for $\reg(G_Q(A))$ in the case $d\ge 2$ is not the best possible. We are interested only in getting a compact bound for $\reg(G_Q(A))$.\medskip

\noindent{\it Remark.} For any $\mm$-primary ideal $I$ of  a generalized Cohen-Macaulay ring $A$ one can bound $\reg(G_I(A))$ by means of $e(I,A)$ and $I(A)$ \cite{Tri} or more general of the extended degree of $I$ \cite{Li}, \cite{RTV1}. Such a bound is not uniform.\medskip

From Theorem \ref{main} we immediately obtain the following uniform bounds for the postulation number. 

\begin{Corollary}
Let $A$ be a generalized Cohen-Macaulay ring with $\dim(A)
= d \geq 1 $. For all parameter ideals $Q$ of $A$,
\begin{align*}
p(Q) & \le \max\{I(A),1\}\ \text{if}\ d = 1,\\
p(Q) & \le \max\{(4I(A))^{(d-1)!}-I(A),1\}\ \text{if}\ d \ge 2.
\end{align*}
\end{Corollary}

Let $R_I(A) = \oplus_{n \ge 0}I^n$ be the Rees algebra of $Q$.
If we represent 
$R_I(A) = A[X_1,...,X_d]/\Im$, where $\Im$ is the set of all forms vanishing at $x_1,...,x_d$, then the maximal degree of the minimal generators of $\Im$ is called the {\it relation type} of $I$, denoted by $\reltype(I)$. It is known \cite[Corollary 1.3 and Proposition 4.1]{Tru} that
$$\reltype(I) \le \reg(R_I(A))+1.$$

On the other hand, Ooishi \cite[Lemma 4.8]{Oo} (see also \cite[Corollary 3.3]{Tru2}) showed that
$$\reg(R_I(A)) = \reg(G_I(A))$$
Therefore, we can deduce from Theorem \ref{main} the following uniform bounds for the relation type. 

\begin{Corollary}
Let $A$ be a generalized Cohen-Macaulay ring with $\dim(A)
= d \geq 1 $. For all parameter ideals $Q$ of $A$,
\begin{align*}
\reltype(Q) & \le \max\{I(A),1\}\ \text{if}\ d = 1,\\
\reltype(Q) & \le \max\{(4I(A))^{(d-1)!}-I(A),1\}\ \text{if}\ d \ge 2.
\end{align*}
\end{Corollary}

Finally, we discuss the problem which local rings have uniform bounds for the relation type of parameter ideals. 
Wang \cite{W1} showed that two-dimensional local rings always have this property. Recently, Aberbach, Ghezzi and Ha \cite{AGH} found out that certain local rings with one-dimensional non-Cohen-Macaulay locus also have  this property. Therefore, the class of local rings  with this property must be larger than the class of generalized Cohen-Macaulay rings. 
However, we shall see that generalized Cohen-Macaulay rings are exactly the local rings for which there is a uniform bound for the relation type of parameter ideals of all quotient rings by ideals generated by subsystems of parameters.

\begin{Theorem}
Let $A$ be a local ring with $d = \dim A \ge 1$. Then the following conditions are equivalent:\par
{\rm (i) } $A$ is a generalized Cohen-Macaulay  ring,\par
{\rm (ii) } There exists an integer $r$ such that for every quotient ring $A/J$, where $J$ is an ideal generated by a subsystem of parameters of $A$, the regularity of the associated rings of all parameters ideals is bounded by $r$,\par
{\rm (ii) } There exists an integer $r$ such that for every quotient ring $A/J$, where $J$ is an ideal generated by a subsystem of parameters of $A$, the relation type of all parameters ideals is bounded by $r$.
\end{Theorem}

\begin{pf}
(i) $\Rightarrow$ (ii). Assume that $A$ is a generalized Cohen-Macaulay ring. Then $A/J$ is also generalized Cohen-Macaulay ring. By Theorem \ref{main}, the regularity of the associated rings of parameters ideals of $A/J$ is bounded above by $(4I(A/J))^{d!} \le (4I(A))^{d!}$.\par

(ii) $\Rightarrow$ (iii). This is a consequence of the formula $\reltype(I) \le \reg(G_I(A))+1$.\par

(iii) $\Rightarrow$ (i). 
Assume that there exists an integer $r$ such that for every quotient ring $B = A/J$, where $J$ is an ideal generated by a subsystem of parameters of $A$, the relation type of parameters ideals of $B$ is bounded by $r$. Consider the case $J = (x_1,...,x_{d-1})$. Let $x = x_d$ be any element of $A$ such that $x_1,...,x_d$ is a system of parameters of $A$. It is easy to check that
$$R_{xB}(B) \cong B[T]/(zT^n|\ z \in 0_B:x^n, n \ge 1).$$
Since $\reltype(xB) \le r$, we must have $0_B:x^n = 0_B:x^r$ for $n > r$. Therefore, 
$x^r(J:x^n) \subseteq J$. Since $ \cup_{n > r}J: x^n =  \cup_{n > r}J: \mm^n$, this implies
$$x^r(\cup_{n \ge r}J: \mm^n)  \subseteq J.$$
Let $s$ denote the minimal number of generators of $\mm$.  Then 
$\mm^{rs}$ is contained in the ideal generated by elements of the form $x^r$. Therefore, $\mm^{rs}(\cup_{n > r}J: \mm^n) \subseteq J$ so that
$$\cup_{n > r}J: \mm^n \subseteq  J:\mm^{rs}.$$
Since $J:x \subseteq \cup_{n > r}J: \mm^n$, we can conclude that
$$(x_1,...,x_{d-1}):x_d \subseteq (x_1,...,x_{d-1}): \mm^{rs}$$
for every system of parameters $x_1,...,x_d$ of $A$. By Lemma \ref{characterization},
this implies the generalized Cohen-Macaulayness of $A$.
\end{pf}

The above theorem shows that if there exists a uniform bound for the relation type of parameter ideals, then the bound must depend on invariants which may increase when we pass to quotient rings by
ideals generated by subsystems of parameters.  It is hard to find such invariants.

\end{document}